\magnification=\magstep1

\centerline 
{A ``COUSIN'' OF A THEOREM OF CS\'AKI AND FISCHER}
\bigskip
\bigskip

\noindent Richard C.\ Bradley \hfil\break
Department of Mathematics \hfil\break
Indiana University \hfil\break
Bloomington \hfil\break
Indiana 47405 \hfil\break
USA \medskip

\noindent bradleyr@indiana.edu
\hfil\break
\bigskip

   {\bf Abstract.}\ \ A 1963 theorem of P.\ Cs\'aki
and J.\ Fischer deals with the ``maximal correlation
coefficient'' in the context of independent pairs of
$\sigma$-fields on a probability space.  
Here a somewhat restricted ``cousin'' of their result
is presented for the same context, but involving in part 
an analogous measure of dependence based only on 
correlations of indicator functions.
It was first proved by the author in an unpublished
1978 Ph.D.\ thesis.
An example is constructed to show a limitation of
this ``cousin''.
Also, this ``cousin'' is used to trivially embellish a 
very sharp 2013 example of R.\ Peyre in connection with 
the comparison of these two measures of dependence. 
\hfil\break
\vfill\eject
   
\noindent {\bf 1.  Introduction}
\bigskip

   Since the papers of Rosenblatt [22] and Ibragimov [16]
and other related works, there has been an extensive development of limit theory under ``strong mixing conditions''.
(For more on such conditions, see e.g.\ [6].)\ \ 
That has motivated a study of ``structural'' properties
of, and connections between, the
strong mixing conditions themselves.
That in turn has motivated a study of the properties of, 
and the connections between, the
various ``measures of dependence'' that form the basis
for such strong mixing conditions.
The ``maximal correlation coefficient'', the measure of dependence which is the basis for the ``$\rho$-mixing
condition'', has been of particular interest. 
Of special interest is a theorem of Cs\'aki and Fischer
[13] involving the maximal correlation coefficient
in the context of independent pairs of $\sigma$-fields.
Here a somewhat restricted ``cousin'' of their result
is presented for the same context, but involving in part 
an analogous measure of dependence based only on 
correlations of indicator functions.
An example is constructed to show a limitation of
this ``cousin''.
Also, this ``cousin'' is used to trivially embellish a 
very sharp example of Peyre [20] in connection with 
the comparison of these two measures of dependence.  
\medskip  

   Suppose $(\Omega, {\cal F}, P)$ is a probability
space.
The indicator function of a given event $A$ will be
denoted by either $I_A$ or $I(A)$.
The term ``$\sigma$-field'' will always refer to a
$\sigma$-field (always $\subset {\cal F}$) on $\Omega$. 
For any two $\sigma$-fields ${\cal A}$ and ${\cal B}$
$\subset {\cal F}$, define the following four 
measures of dependence:
First, 
$$  
\psi({\cal A}, {\cal B})\ =\
\sup_{A \in {\cal A}, \thinspace B \in {\cal B}}
{{|P(A \cap B) - P(A)P(B)|} \over {P(A)P(B)}}\ .  
\eqno (1.1)
$$
Next,
$$\lambda({\cal A}, {\cal B})\ :=\
\sup_{A \in {\cal A}, \thinspace B \in {\cal B}}
{{|P(A \cap B) - P(A)P(B)|} \over 
{[P(A)P(B)]^{1/2}}}\ .  
\eqno (1.2) 
$$
Next,
$$ \eqalignno{
\tau({\cal A}, {\cal B})\ &:=\
\sup_{A \in {\cal A}, \thinspace B \in {\cal B}}
|{\rm Corr}(I_A, I_B)|\ \cr
 &\ =\ 
\sup_{A \in {\cal A}, \thinspace B \in {\cal B}}
{{|P(A \cap B) - P(A)P(B)|} \over 
{[P(A) \cdot (1 - P(A)) \cdot 
P(B) \cdot (1 - P(B))]^{1/2}}}\ .  
& (1.3) \cr 
}$$
Finally, 
$$ \rho({\cal A}, {\cal B})\ :=\ \sup |{\rm Corr}(X,Y)|
\eqno (1.4) $$
where the supremum in (1.4) is taken over all pairs of
square-integrable random variables $X$ and $Y$ such that
$X$ is ${\cal A}$-measurable and    
$Y$ is ${\cal B}$-measurable.
In (1.1), (1.2), and (1.3), the fraction $0/0$ is 
interpreted as 0.
In (1.3), the second equality is a standard elementary
calculation.
In (1.3) and (1.4), ``Corr'' denotes the correlation.
\medskip

   The quantity $\psi({\cal A}, {\cal B})$ 
was implicitly present in work of 
Doeblin [14] 
involving a ``continued fraction process''; 
for some details, see e.g.\ Iosifescu [17].
Later on, more explicitly, that quantity 
$\psi({\cal A}, {\cal B})$ was, for general stochastic
processes, the basis for the
*-mixing condition in Blum, Hanson, and
Koopmans [1] and for the $\psi$-mixing
condition in Philipp [21] and other papers.       
\medskip

   The quantity $\rho({\cal A}, {\cal B})$ in (1.4) is the 
well known ``maximal correlation'' coefficient, first 
studied by Hirschfeld [15].
It was, for stochastic processes, the basis for the 
$\rho$-mixing condition, 
introduced by Kolmogorov and Rozanov [18].
\medskip

   The measures of dependence 
$\lambda({\cal A}, {\cal B})$ and
$\tau({\cal A}, {\cal B})$, formulated exactly as in
(1.2) and (1.3), were examined in [3], [9], and [12]
with a view toward allowing arguments involving the
maximal correlation coefficient $\rho({\cal A}, {\cal B})$ 
to be simplified by the converting of such arguments 
from pairs of 
(square-integrable) random variables to pairs of events.
With that in mind, let us look at the comparison of these 
three measures of dependence.
First, the following inequalities hold:
$$
\lambda({\cal A}, {\cal B})\ 
\leq\ \tau({\cal A}, {\cal B})\
\leq\ \rho({\cal A}, {\cal B})\
\leq\ 1 \quad \bigl( {\rm and\ also} \quad 
\rho({\cal A}, {\cal B})\ \leq 
\psi({\cal A}, {\cal B})\ \bigl) . \eqno (1.5)    
$$       
The first three are trivial; and the last one is well
known and elementary (see e.g.\ 
[6, v1, Proposition 3.11(b)]).
By a simple calculation, for any two events $A$ and $B$,
the quantity $|P(A \cap B) - P(A)P(B)|$ remains 
unchanged if either $A$ or $B$ is replaced by its 
complement.
Consequently, the definition in (1.2) 
(as well as those in (1.1) and (1.3)) does not change if
one restricts to pairs of events $A \in {\cal A}$ 
and $B \in {\cal B}$ such that $P(A) \leq 1/2$ and
$P(B) \leq 1/2$.
As a simple consequence, one also has the following inequality:
$$ \tau({\cal A}, {\cal B})\ \leq\ 
2 \cdot \lambda({\cal A}, {\cal B})\ . 
\eqno (1.6) $$
The author [3] proved the crude inequality
$\rho({\cal A}, {\cal B}) \leq 
13 \cdot [\tau({\cal A}, {\cal B})]^{1/31}$.
Together with (1.6) and the first two inequalities 
in (1.5), that showed that the three measures of dependence
$\lambda(.\, ,.)$, $\tau(.\, ,.)$, and $\rho(.\, ,.)$ are
``equivalent'', in that they all become arbitrarily small
as any one of them becomes sufficiently small.
Later, the author and Bryc [9, Theorem 1.1(ii)], and independently Bulinskii [12, the Theorem], 
showed that there exists a
universal positive constant $C$ such that the inequality
$$ \rho({\cal A}, {\cal B})\ \leq\ 
C \cdot \lambda({\cal A}, {\cal B}) \cdot
[1 - \log\thinspace \lambda({\cal A}, {\cal B})]
\eqno (1.7) $$
always holds.
For a quite gentle proof of that result, adapted partly
from Bulinskii's [12] very sharp improvement
of the crude calculations in [3], 
see [6, v1, Theorem 4.15].
The author, Bryc, and Janson [10, Theorem 3.1] 
showed (as a special case of a more general result)
that the inequality in (1.7) is within a constant factor of
being sharp --- i.e.\ that there 
exists a universal positive
constant $A$ such that the following holds:
For any $t \in [0, 1]$, there exist a probability 
space $(\Omega, {\cal F}, P)$
and $\sigma$-fields 
${\cal A}$ and ${\cal B}\ \subset {\cal F}$
such that $\tau({\cal A}, {\cal B}) \leq t$ and
$\rho({\cal A}, {\cal B}) \geq 
A \cdot t \cdot (1 - \log t)$.
For a more gentle proof of that particular result, see  
[6, v1, Theorem 4.16].
More recently, with a much more sophisticated argument,
an ``exact'' (i.e.\ ``best possible'') version of (1.7)
(involving the measure of dependence $\tau(. \, , .)$) 
was proved by Peyre [20]; that result will be stated in 
Theorem 2 below.          
\medskip

   (In connection with (1.7), note that by simple 
calculus, the expression $t(1 - \log t)$ is strictly
increasing as $t$ increases in $[0,1]$. 
Here and below, $0 \log 0 := 0$.)
\medskip
   
   Via inequalities such as (1.7), the measures of 
dependence $\lambda(.\, ,.)$ and 
$\tau(.\, ,.)$ are useful in simplifying some
arguments pertaining to the $\rho$-mixing
condition and other conditions based on the
maximal correlation coefficient $\rho(.\, ,.)$.
The measure of dependence $\lambda(.\, ,.)$ is of
course the easiest of the three to work with; and it has
been used by the author [4][5][7][8] 
to simplify proofs of the following results: 
(1) the equivalence of the Rosenblatt [22] ``strong
mixing condition'' with a certain condition of 
``$\rho$-mixing except on small sets''
(a phrase coined by Magda Peligrad, who had originally 
brought that latter condition to the author's attention);
(2) the ``$\rho^*$-mixing'' property 
(the stronger, ``interlaced'' variant of $\rho$-mixing)
for certain Markov chains, including as a
special case the strictly stationary, finite-state, irreducible, aperiodic Markov chains, 
(3) the $\rho^*$-mixing property of INAR
(``integer-valued autoregressive'') processes of order 1
with ``Poisson innovations''; and
(4) the existence of strictly stationary,
countable-state, reversible Markov chains that 
satisfy $\rho$-mixing
(and hence also geometric ergodicity) but fail to
satisfy $\rho^*$-mixing.
(For more on (1) and (2), see also
[6, v2, pp.\ 415-423, and v1, Theorem 7.15].)
\medskip

   Also, when Magda Peligrad formulated, and
developed some central limit theory under,
a ``two-part'' mixing condition --- in essence a 
``hybrid'' of the (Rosenblatt [22]) strong mixing 
condition and the $\rho$-mixing condition ---
she adapted the measure of dependence $\lambda(.\, ,.)$
to simplify the formulation of the ``component'' of her 
two-part mixing condition that was 
related to $\rho$-mixing.
For details, see [11] and [19] and also
[6, v2, Chapter 18].
\medskip

   Recall from (1.5) and (1.6) that the measures of
dependence $\lambda(.\, ,.)$ and $\tau(.\, ,.)$
differ from each other by a most a factor of 2
--- and hence trivially (1.7) holds (with at most
a change in the constant factor $C$) with
$\lambda(.\, ,.)$ replaced by $\tau(.\, ,.)$.
Of these two measures of dependence, the latter one
seems better suited for making ``exact comparisons''
with the maximal correlation coefficient $\rho(.\, ,.)$.
Trivially, for any two $\sigma$-fields 
${\cal A}$ and ${\cal B}$ that are each purely
atomic with exactly two atoms, the equality
$\rho({\cal A}, {\cal B}) = 
\tau({\cal A}, {\cal B})$ holds.   
With a quite elementary argument, 
in the case where one of the 
$\sigma$-fields is purely atomic with exactly two atoms 
and the other $\sigma$-field is ``unrestricted'',
the author and Bryc [9, Theorem 4.3 and Example 4.4]
derived the following result, giving 
an ``exact'' (i.e.\ ``best possible'') inequality:
\hfil\break

   {\bf Theorem 1} ([9, Theorem 4.3 and Example 4.4]).\ \ 
(I)\ {\sl Suppose
$(\Omega, {\cal F}, P)$ is a probability space,
${\cal A}$ and ${\cal B}$ are 
$\sigma$-fields $\subset {\cal F}$, and
the $\sigma$-field ${\cal A}$ is purely atomic with
exactly two atoms; then}
$$ \rho({\cal A}, {\cal B})\ \leq\ 
\tau({\cal A}, {\cal B}) \cdot
[1 - \log\thinspace \tau({\cal A}, {\cal B})]^{1/2}\ .
\eqno (1.8) $$ 
(II)\ {\sl For any $t \in [0,1]$ and any $a \in (0,1)$,
there exist a probability space $(\Omega, {\cal F}, P)$
and $\sigma$-fields 
${\cal A}$ and ${\cal B}\ \subset {\cal F}$
such that 
(i) ${\cal A}$ is purely atomic with exactly two
atoms $A$ and $A^c$, such that 
$P(A) = a$ and $P(A^c) = 1-a$, and
(ii) $\tau({\cal A}, {\cal B}) = t$ and
$\rho({\cal A}, {\cal B}) = 
t \cdot (1 - \log t)^{1/2}$.}
\bigskip

    The ``sharp constant'' in (1.8) is of course 
(implicitly) 1.
Note that in the inequality in (1.8), the ``log term''
has an exponent $1/2$ that is not present in (1.7).
That exponent in (1.8) is of course connected with the
extra restriction (not present in (1.7))
that one of the $\sigma$-fields is 
purely atomic with exactly two atoms.
\medskip

   (Again, by simple calculus, the expression
$t(1 - \log t)^{1/2}$ is strictly increasing as $t$
increases in $[0,1]$.)
\medskip 

   In the original context involving no restriction on 
either $\sigma$-field, 
Peyre [20, Theorem 3.1 and Theorem 4.1] showed with 
a much more sophisticated argument that in the version of (1.7) with $\lambda(.\, ,.)$ replaced by $\tau(.\, ,.)$, 
the ``sharp constant'' is again 1 and the
resulting inequality is ``exact'' (i.e.\ ``best possible'').
Here we shall state his result in the notations 
used here in this paper.
(The notations used by Peyre [20] slightly conflict 
with those used here.)
\hfil\break

{\bf Theorem 2} ([20, Theorems 3.1 and 4.1]). \ \ (I)\ {\sl Suppose
$(\Omega, {\cal F}, P)$ is a probability space, and
${\cal A}$ and ${\cal B}$ are 
$\sigma$-fields $\subset {\cal F}$; then}
$$ \rho({\cal A}, {\cal B})\ \leq\ 
\tau({\cal A}, {\cal B}) \cdot
[1 - \log\thinspace \tau({\cal A}, {\cal B})]\ .
\eqno (1.9) $$ 
(II)\ {\sl For any $t \in (0,1)$ and any 
$\varepsilon \in (0,t)$,
there exist a probability space $(\Omega, {\cal F}, P)$
and $\sigma$-fields 
${\cal A}$ and ${\cal B}\ \subset {\cal F}$
such that 
$\tau({\cal A}, {\cal B}) \leq t$ and
$\rho({\cal A}, {\cal B}) > 
[t \cdot (1 - \log t)] - \varepsilon$.}
\bigskip

   Obviously (1.9) is a ``sharpest possible''
version of (1.7)
(with $\lambda(.\, ,.)$ replaced by $\tau(.\, ,.)$); 
and the example described in (II)
here gives a very sharp improvement compared to the 
special case of the example in [10] that was
alluded to right after (1.7).
\medskip

   In Corollary 5 below, we shall show that
a variant or ``cousin'' of a
result of Cs\'aki and Fischer [13, Theorem 6.2] 
allows one to trivially ``embellish'' Peyre's example described in Theorem 2(II) in such a way that (also) 
there exist events
$A \in {\cal A}$ and $B \in {\cal B}$ such that
${\rm Corr}(I_A,I_B) = t$ (recall (1.3)).  
So far, apparently no way has been found to also
achieve ``equality in (1.9)'' in such an example,
that is, to achieve the equality
$\rho({\cal A}, {\cal B}) = [t \cdot (1 - \log t)]$.
\medskip

   First let us state the result of Cs\'aki and
Fischer [13, Theorem 6.2] itself: 
\bigskip

   {\bf Theorem 3} ([13, Theorem 6.2]).\ \ 
{\sl Suppose $(\Omega, {\cal F}, P)$ is a probability space,
and 
${\cal A}_n$ and ${\cal B}_n$,
$n \in {\bf N}$ are $\sigma$-fields 
$\subset {\cal F}$
such that the $\sigma$-fields 
${\cal A}_n \vee {\cal B}_n$, 
$n \in {\bf N}$ are independent.
Then 
$$ \rho \biggl(\ \bigvee_{n \in {\bf N}}{\cal A}_n,\ 
\bigvee_{n \in {\bf N}}{\cal B}_n \biggl)\ =\ 
\sup_{n \in {\bf N}}\rho({\cal A}_n, {\cal B}_n)\ .
\eqno (1.10) 
$$}

   For a generously detailed proof of this theorem
(essentially, an induction argument given by Witsenhausen [23], followed by a standard measure-theoretic argument,
all with plenty of detail), see [6, v1, Theorem 6.1].
Theorem 3 has been used in the proofs of results
in [5][7][8][13][23] and many other papers, as well as 
in the proofs of numerous results in [6].  
\medskip   

   Now let us look at a variant or ``cousin'' 
--- in some limited sense ---
of Theorem 3.
The following result was stated and proved years ago 
by the author
(in an equivalent form, without explicit use of the 
notations $\tau(.\, ,.)$ and $\psi(.\, ,.)$) 
in [2, Theorem 6], in an unpublished Ph.D.\ 
thesis.
\hfil\break

   {\bf Theorem 4} ([2, Theorem 6]).\ \ {\sl Suppose  
$(\Omega, {\cal F}, P)$ is a probability space,
and 
${\cal A}_1$, ${\cal B}_1$, ${\cal A}_2$, and
${\cal B}_2$ are $\sigma$-fields $\subset {\cal F}$
such that 
the $\sigma$-fields 
${\cal A}_1 \vee {\cal B}_1$ and 
${\cal A}_2 \vee {\cal B}_2$ are independent.
Then
$$ \tau({\cal A}_1 \vee {\cal A}_2, \thinspace 
{\cal B}_1 \vee {\cal B}_2)\ \leq\ 
\max\{ \tau({\cal A}_1, {\cal B}_1), \thinspace
\psi({\cal A}_2, {\cal B}_2) \}\ . 
\eqno (1.11) $$}

   The proof (from [2]) of Theorem 4 will be given
in Section 2.
A limitation of Theorem 4 in connection with the term
$\psi({\cal A}_2, {\cal B}_2)$ in (1.11) will be treated 
in Theorem 6 below.
\medskip

   Of course by Theorem 4 and induction, followed by a 
standard measure-theoretic argument, 
as an analog of (1.10), one has that
under the hypothesis of Theorem 3, 
$$ \tau \biggl(\ \bigvee_{n \in {\bf N}}{\cal A}_n,\ 
\bigvee_{n \in {\bf N}}{\cal B}_n \biggl)\ \leq\ 
\sup \Bigl\{ \tau({\cal A}_1, {\cal B}_1),\  
\sup_{n \geq 2} \psi({\cal A}_n, {\cal B}_n)\Bigl\}\ .
\eqno (1.12)  
$$ 

   As an application of Theorem 4, the example given
by Peyre [20, Theorem 4.1] described in Theorem 2(II) 
will be trivially ``embellished'', in a certain way 
alluded to above:   
\hfil\break

   {\bf Corollary 5} (trivial embellishment of Peyre's 
[20] example).
{\sl For any $t \in (0,1)$ and any 
$\varepsilon > 0$,
there exist a probability space $(\Omega, {\cal F}, P)$
and $\sigma$-fields 
${\cal A}$ and ${\cal B}\ \subset {\cal F}$
with the following properties:
\hfil\break
(i) there exist events $A \in {\cal A}$ and
$B \in {\cal B}$ such that 
$\tau({\cal A}, {\cal B}) = {\rm Corr}(I_A,I_B) = t$; and
\hfil\break
(ii) $\rho({\cal A}, {\cal B}) > 
[t \cdot (1 - \log t)] - \varepsilon$.}
\bigskip

   The proof of Corollary 5 will be given in 
Section 3 below.
Its proof will make critical use of Peyre's [20]
example itself (as described in Theorem 2(II)).
\medskip

   The final result, Theorem 6 below, will show that
in Theorem 4, in eq.\ (1.11), the term 
$\psi({\cal A}_2, {\cal B}_2)$ cannot be replaced by
$\tau({\cal A}_2, {\cal B}_2)$.
It seems to be an open question whether or not in 
Theorem 4, in eq.\ (1.11), that term 
$\psi({\cal A}_2, {\cal B}_2)$ can be replaced by
$\rho({\cal A}_2, {\cal B}_2)$.
\hfil\break

   {\bf Theorem 6.}\ \ {\sl Suppose
$$ 0 < t < 1.  \eqno (1.13) $$
Then there exist a probability space
$(\Omega, {\cal F}, P)$ and 
$\sigma$-fields ${\cal A}_1$, ${\cal B}_1$,
${\cal A}_2$, and ${\cal B}_2$ ($\subset {\cal F}$)
with the following properties:  
$${\cal A}_1 \vee {\cal B}_1\ {\rm and}\  
{\cal A}_2 \vee {\cal B}_2\ {\rm are\ independent}
\eqno (1.14) $$
and
$$ \tau({\cal A}_1, {\cal B}_1)\ =\ 
\tau({\cal A}_2, {\cal B}_2)\ =\ t\ , \eqno (1.15) $$
but 
$$ \tau({\cal A}_1 \vee {\cal A}_2, \thinspace 
{\cal B}_1 \vee {\cal B}_2)\ >\ t\ . \eqno (1.16) $$}

   Theorem 6 will be proved in Section 4.
Its proof will make critical use of Peyre's [20] example.
(For sufficiently small $t$, the proof could instead
use in a similar way the example from [10] alluded to
right after (1.7).) 
\hfil\break

\noindent {\bf 2. Proof of Theorem 4}
\bigskip

The proof is essentially as given by the 
author [2, Theorem 6 (its proof)].
The proof will be divided into several small ``steps''.
\hfil\break 
   
   {\bf Step 1.}\ \ As in the statement of 
Theorem 4, suppose
$(\Omega, {\cal F}, P)$ is a probability space,
and 
${\cal A}_1$, ${\cal B}_1$, ${\cal A}_2$, and
${\cal B}_2$ are $\sigma$-fields ($\subset {\cal F}$)
such that  
$${\cal A}_1 \vee {\cal B}_1\ {\rm and}\  
{\cal A}_2 \vee {\cal B}_2\ {\rm are\ independent}.
\eqno (2.1) $$

If $\psi({\cal A}_2, {\cal B}_2) = \infty$, then 
(1.11) holds trivially and we are done.
Therefore we assume that 
$\psi({\cal A}_2, {\cal B}_2) < \infty$.
With a reminder of that assumption built in 
(see also (1.5)),  
define the nonnegative quantity
$$ \theta\ :=\ \max\{ \tau({\cal A}_1, {\cal B}_1), 
\thinspace
\psi({\cal A}_2, {\cal B}_2) \}\ 
<\ \infty\ . \eqno (2.2) $$
Our task is to show that
$\tau({\cal A}_1 \vee {\cal A}_2, \thinspace 
{\cal B}_1 \vee {\cal B}_2) \leq \theta$.
\medskip

   Refer to (1.3). Suppose 
$$A_0 \in {\cal A}_1 \vee {\cal A}_2 \quad {\rm and} \quad
B_0 \in {\cal B}_1 \vee {\cal B}_2\ .  \eqno (2.3) $$
It suffices to prove that 
$$ |P(A_0 \cap B_0) - P(A_0)P(B_0)|\ \leq\ 
\theta \cdot 
[P(A_0) \cdot (1 - P(A_0)) \cdot P(B_0) 
\cdot (1 - P(B_0))]^{1/2}\ .$$
Suppose $\varepsilon > 0$.  
It suffices to prove that 
$$ |P(A_0 \cap B_0) - P(A_0)P(B_0)|\ \leq\ 
\theta \cdot 
[P(A_0) \cdot (1 - P(A_0)) \cdot P(B_0) 
\cdot (1 - P(B_0))]^{1/2}\ +\varepsilon\ . \eqno (2.4) $$  

   {\bf Step 2.}\ \ We shall first make a long statement
(ending with eq.\ (2.6) below), 
and then briefly justify it.
By (2.3) and a standard measure-theoretic argument, 
there exist
events $A^* \in {\cal A}_1 \vee {\cal A}_2$ and
$B^* \in {\cal B}_1 \vee {\cal B}_2$ with the following
three properties (P1), (P2), and (P3):
\medskip

\noindent (P1) $A^* = \bigcup_{i=1}^I(C_i \cap D_i)$ where 
\hfil\break
(i) $I$ is a positive integer, \hfil\break
(ii) $C_i \in {\cal A}_1$ and $D_i \in {\cal A}_2$
for each $i \in \{1, 2, \dots, I\}$, and \hfil\break
(iii) the events $D_1, D_2, \dots, D_I$ together
form a partition of the sample space $\Omega$.
\medskip

\noindent (P2) $B^* = \bigcup_{j=1}^J(E_j \cap F_j)$ where 
\hfil\break
(i) $J$ is a positive integer, \hfil\break
(ii) $E_j \in {\cal B}_1$ and $F_j \in {\cal B}_2$
for each $j \in \{1, 2, \dots, J\}$, and \hfil\break
(iii) the events $F_1, F_2, \dots, F_J$ together
form a partition of the sample space $\Omega$.
\medskip

\noindent (P3) One has that
$$ \Bigl|\ 
|P(A^* \cap B^*) - P(A^*)P(B^*)| 
- |P(A_0 \cap B_0) - P(A_0)P(B_0)|\ \Bigl|\  
\leq\ \varepsilon/2\   
\eqno (2.5) $$
and (recall the ``$< \infty$'' in (2.2))
$$ \eqalignno{
\Bigl|\ 
\theta \cdot 
[P(A^*) \cdot &(1 - P(A^*)) \cdot P(B^*) 
\cdot (1 - P(B^*))]^{1/2}\ \cr
-\ 
& \theta \cdot 
[P(A_0) \cdot (1 - P(A_0)) \cdot P(B_0) 
\cdot (1 - P(B_0))]^{1/2}\
\Bigl|\ 
\leq\ \varepsilon/2\ . & (2.6) \cr
}$$

   Let us briefly review the justification
of this assertion involving properties 
(P1), (P2), and (P3):
\medskip

   Let the symmetric difference of any two events 
$G$ and $H$ be denoted $G \triangle H$.
It is well known that (i) for any two events $G$ and $H$,
one has that $|P(G)-P(H)| \leq P(G \triangle H)$, and 
(ii) for any four events
$G_1$, $H_1$, $G_2$, and $H_2$, one has that
$P((G_1 \cap H_1) \triangle (G_2 \cap H_2))
\leq P(G_1 \triangle G_2) + P(H_1 \triangle H_2)$.
\medskip

   By the first part of (2.3) and a well known
measure theoretic argument, for each 
$\gamma > 0$, there exist positive
integers $K$ and $I$ and partitions 
$\{G_1, G_2, \dots, G_K\}$ and 
$\{D_1, D_2, \dots, D_I\}$ of $\Omega$, with
$G_k \in {\cal A}_1$ for each $k \in \{1, 2, \dots, K\}$
and $D_i \in {\cal A}_2$ for each 
$i \in \{1, 2, \dots, I\}$, and an event
$A^*$ which is the union of some (or all or none) of the ``rectangles'' $G_k \cap D_i$, such that
$P(A^* \triangle A_0) \leq \gamma$.
Then for each fixed $i \in \{1, 2, \dots, I\}$, one can let
$C_i := \bigcup_k G_k$ where the union 
(possibly empty) is taken
over all $k \in \{1, 2, \dots, K\}$ such that
$G_k \cap D_i \subset A^*$.
Then the set $A^*$ has the form specified in property (P1)
(and satisfies $P(A^* \triangle A_0) \leq \gamma$).
Similarly for each $\gamma > 0$, there exists an event
$B^*$ satisfying property (P2) such that
$P(B^* \triangle B_0) \leq \gamma$.
By taking $\gamma > 0$ sufficiently small, and using
observations (i) and (ii) in the preceding paragraph, 
one can ensure that property (P3) 
(both eqs.\ (2.5) and (2.6) --- recall again the 
``$< \infty$'' in (2.2)) holds as well.
\medskip

   To prove (2.4) and thereby complete the proof of
Theorem 4, it now suffices to prove that the events
$A^*$ and $B^*$, satisfying properties (P1), (P2), and
(P3), satisfy
$$ |P(A^* \cap B^*) - P(A^*)P(B^*)|\ \leq\ 
\theta \cdot 
[P(A^*) \cdot (1 - P(A^*)) \cdot P(B^*) 
\cdot (1 - P(B^*))]^{1/2}\ .    \eqno (2.7) $$

   {\bf Step 3.}\ \ By property (P1)(iii) (in Step 2),
the events $C_i \cap D_i$, for different values of $i$,
are (pairwise) disjoint.
Thus from property (P1), the events 
$C_i  \cap D_i$, $i \in \{1, 2, \dots, I\}$ 
(some of those events may be empty) form a
partition of the event $A^*$.
Similarly from property (P2), the events
$E_j \cap F_j$, $j \in \{1, 2, \dots, J\}$ 
(some of those events may be empty) form a
partition of the event $B^*$.
It follows that the events 
$C_i \cap D_i \cap E_j \cap F_j$,
$(i,j) \in \{1, 2, \dots, I\} \times \{1, 2, \dots, J\}$
(some of those events may be empty)
form a partition of the event $A^* \cap B^*$.
It follows that
$$ \eqalignno{
P&(A^* \cap B^*) - P(A^*)P(B^*) \cr 
&=\ 
\sum_{i=1}^I \sum_{j=1}^J P(C_i \cap D_i \cap E_j \cap F_j)
\ -\ 
\biggl[\thinspace \sum_{i=1}^I P(C_i \cap D_i)\biggl] \cdot 
\biggl[\thinspace \sum_{j=1}^J P(E_j \cap F_j)\biggl]\ . 
\indent & (2.8) \cr
}$$
Applying properties (P1)(ii) and (P2)(ii) and eq.\ (2.1)
to (2.8), one obtains  
$$ \eqalignno{
P&(A^* \cap B^*) - P(A^*)P(B^*) \cr
&=\ 
\sum_{i=1}^I \sum_{j=1}^J P(C_i \cap E_j)P(D_i \cap F_j)
\ -\ 
\thinspace \sum_{i=1}^I \thinspace \sum_{j=1}^J
P(C_i)P(D_i)P(E_j)P(F_j) \cr
&=\
\sum_{i=1}^I \sum_{j=1}^J P(D_i \cap F_j) \cdot
[P(C_i \cap E_j) - P(C_i)P(E_j)] \cr
& \indent \indent +\ 
\sum_{i=1}^I \sum_{j=1}^J [P(D_i \cap F_j) - P(D_i)P(F_j)]
\cdot P(C_i)P(E_j)\ . & (2.9) \cr  
}$$

   {\bf Step 4.}\ \ The next task is to obtain a useful
alternative formulation of the very last double sum in
(2.9).
For that purpose, let us make some observations.
First, by properties (P1)(ii) and (P2)(ii) and eq.\ (2.1),
followed by property (P2)(iii),  
$$ \eqalignno{ \sum_{i=1}^I \sum_{j=1}^J &
[P(D_i \cap F_j) - P(D_i)P(F_j)] 
\cdot P(C_i) \cdot P(B^*) \cr 
&= P(B^*) \cdot \sum_{i=1}^I \sum_{j=1}^J
[P(C_i \cap D_i \cap F_j) - P(C_i \cap D_i)P(F_j)] \cr
&= P(B^*) \cdot \sum_{i=1}^I 
[P(C_i \cap D_i) - P(C_i \cap D_i)]\ =\ 0\ . 
& (2.10) \cr  
}$$   
Next, by an exactly analogous argument, this time
finishing with an application of property (P1)(iii),      
$$ \eqalignno{ \sum_{i=1}^I \sum_{j=1}^J &
[P(D_i \cap F_j) - P(D_i)P(F_j)] 
\cdot P(A^*) \cdot P(E_j) \cr
&= P(A^*) \cdot \sum_{i=1}^I \sum_{j=1}^J
[P(D_i \cap E_j \cap F_j) - P(D_i)P(E_j \cap F_j)] \cr
&= P(A^*) \cdot \sum_{j=1}^J 
[P(E_j \cap F_j) - P(E_j \cap F_j)]\ =\ 0\ . 
& (2.11) \cr  
}$$
Also of course by properties (P1)(iii) and (P2)(iii),
$$ \eqalignno{ 
\sum_{i=1}^I \sum_{j=1}^J & 
[P(D_i \cap F_j) - P(D_i)P(F_j)] \cdot P(A^*)P(B^*) \cr 
&=\ P(A^*)P(B^*) \cdot \biggl[
\sum_{i=1}^I \sum_{j=1}^J 
P(D_i \cap F_j) 
-\ \sum_{i=1}^I \sum_{j=1}^J 
P(D_i)P(F_j) \biggl] \cr 
&=\ P(A^*)P(B^*) \cdot [1\ -\ 1]\ =\ 0\ . & (2.12) \cr
}$$
Now by incorporating (2.10), (2.11), and (2.12) into
the very last double sum in (2.9), one obtains
from (2.9) itself that
$$ \eqalignno{
P&(A^* \cap B^*) - P(A^*)P(B^*) \cr
&=\
\sum_{i=1}^I \sum_{j=1}^J P(D_i \cap F_j) \cdot
[P(C_i \cap E_j) - P(C_i)P(E_j)] \cr
& \indent +\ 
\sum_{i=1}^I \sum_{j=1}^J [P(D_i \cap F_j) - P(D_i)P(F_j)]
\cdot [P(C_i) - P(A^*)] \cdot [P(E_j) - P(B^*)]\ .   \cr  
}$$
Hence by the triangle inequality and then (2.2), (1.3),
and (1.1) (and properties (P1)(ii) and (P2)(ii)), 
$$ \eqalignno{
|P&(A^* \cap B^*) - P(A^*)P(B^*)| \cr
&\leq\
\sum_{i=1}^I \sum_{j=1}^J P(D_i \cap F_j) \cdot
|P(C_i \cap E_j) - P(C_i)P(E_j)| \cr
& \indent +\ 
\sum_{i=1}^I \sum_{j=1}^J |P(D_i \cap F_j) - P(D_i)P(F_j)|
\cdot |P(C_i) - P(A^*)| \cdot |P(E_j) - P(B^*)| \cr   
& \leq\ 
\sum_{i=1}^I \sum_{j=1}^J P(D_i \cap F_j) \cdot
\theta \cdot [P(C_i) \cdot (1 - P(C_i)) \cdot 
P(E_j) \cdot (1 - P(E_j))]^{1/2} \cr
& \indent +\ 
\sum_{i=1}^I \sum_{j=1}^J \theta \cdot P(D_i)P(F_j)
\cdot |P(C_i) - P(A^*)| \cdot |P(E_j) - P(B^*)|\ . 
 & (2.13) \cr  
}$$ 

   {\bf Step 5.}\ \ Now we shall apply to (2.13) the
Cauchy-Schwarz Inequality.
To put this a little informally, think of a discrete 
measure space with exactly $2IJ$
points, with a positive measure that assigns masses
$\theta \cdot P(D_i \cap F_j)$ respectively to the 
``first $IJ$ points'' and masses 
$\theta \cdot P(D_i)P(F_j)$ to the ``other $IJ$ points''.
With that interpretation (in an obvious form), 
applying the Cauchy-Schwarz
Inequality to (2.13), and then applying
properties (P2)(iii) and (P1)(iii), one obtains
 
$$ \eqalignno{
|P&(A^* \cap B^*) - P(A^*)P(B^*)| \cr  
& \leq\ 
\Biggl[ \sum_{i=1}^I \sum_{j=1}^J 
\theta \cdot P(D_i \cap F_j) 
\cdot [P(C_i) \cdot (1 - P(C_i))] \cr  
& \indent \indent \indent +\ 
\sum_{i=1}^I \sum_{j=1}^J \theta \cdot P(D_i)P(F_j)
\cdot [P(C_i) - P(A^*)]^2 \Biggl]^{1/2} 
\cr
& \indent \cdot\ 
\Biggl[ \sum_{i=1}^I \sum_{j=1}^J 
\theta \cdot P(D_i \cap F_j) \cdot  
[P(E_j) \cdot (1 - P(E_j))] \cr
& \indent \indent \indent \indent +\ 
\sum_{i=1}^I \sum_{j=1}^J \theta \cdot P(D_i)P(F_j)
 \cdot [P(E_j) - P(B^*)]^2 \Biggl]^{1/2} 
\cr
&\ =\ \theta \cdot  
\Biggl[ \sum_{i=1}^I  
 P(D_i) 
\cdot [P(C_i) \cdot (1 - P(C_i))] \ +\ 
\sum_{i=1}^I  P(D_i)
\cdot [P(C_i) - P(A^*)]^2 \Biggl]^{1/2} 
\cr
& \indent \quad \cdot\ 
\Biggl[  \sum_{j=1}^J 
 P(F_j) \cdot  
[P(E_j) \cdot (1 - P(E_j))] \ +\ 
\sum_{j=1}^J  P(F_j)
\cdot [P(E_j) - P(B^*)]^2 \Biggl]^{1/2}  
\cr
&\ =\ \theta \cdot  
\Biggl[ \sum_{i=1}^I  
 P(D_i) 
\cdot [P(C_i) -2 P(C_i)P(A^*) + (P(A^*))^2] 
\Biggl]^{1/2}  
\cr
& \indent \quad \cdot\ 
\Biggl[  \sum_{j=1}^J 
 P(F_j) \cdot  
[P(E_j) -2 P(E_j)P(B^*) + (P(B^*))^2]  
\Biggl]^{1/2}\ . & (2.14)   
\cr     
}$$

   {\bf Step 6.}\ \ Now let us look at the very
last product in (2.14).
To start off, note that by (2.1) and 
property (P1) (see both (ii) and (iii) there)
$$ \eqalignno{
\sum_{i=1}^I & P(D_i) 
\cdot [P(C_i) -2 P(C_i)P(A^*) + (P(A^*))^2] \cr
&=\ (1 - 2P(A^*))\sum_{i=1}^I P(D_i)P(C_i)\ 
+\ (P(A^*))^2 \sum_{i=1}^I  P(D_i) \cr
\cr 
&=\ (1 - 2P(A^*))\sum_{i=1}^I P(D_i \cap C_i)\ 
+\ (P(A^*))^2 \cdot 1 \cr
&= (1 - 2P(A^*))\cdot P(A^*)\ 
+\ (P(A^*))^2\ =\ P(A^*) \cdot (1 - P(A^*))\ . 
& (2.15) \cr         
}$$  
By an exactly analogous argument, using 
property (P2) instead of (P1), one has that
$$\sum_{j=1}^J P(F_j) \cdot  
[P(E_j) -2 P(E_j)P(B^*) + (P(B^*))^2]\ =\ 
P(B^*) \cdot (1 - P(B^*))\ . \eqno (2.16) $$
Applying (2.15) and (2.16) to the last product in (2.14),
one obtains from (2.14) itself that (2.7) holds.
That completes the proof of Theorem 4.
\hfil\break

\noindent {\bf 3. Proof of Corollary 5}
\bigskip

   As in the statement of Corollary 5, 
suppose $t \in (0,1)$ and $\varepsilon > 0$.
\medskip 
      
   The construction given by Peyre [20, Theorem 4.1],
described in Theorem 2(II), can be interpreted in the following way:  
On some probability space $(\Omega, {\cal F}, P)$,
there exists a random vector $(X_1, X_2)$ such that
$$ \eqalignno{ 
\tau(\sigma(X_1), \sigma(X_2))\ & \leq\ t 
\quad {\rm and} & (3.1) \cr
\rho(\sigma(X_1), \sigma(X_2))\ & >\  
[t \cdot (1 - \log t)] - \varepsilon\ . & (3.2) \cr
}$$      
Here and below, the notation $\sigma(\dots)$ means the
$\sigma$-field generated by $(\dots)$.
The details of Peyre's construction, which is
intricate and quite long, need not be spelled
out here.
\medskip

   Enlarging the probability space if necessary, let
$(Y_1, Y_2)$ be a random vector which is independent of
the random vector $(X_1, X_2)$ and has the following
distribution:  For each element
$(y_1, y_2) \in \{-1,1\} \times \{-1,1\}$,  
$$ P \bigl((Y_1, Y_2) = (y_1, y_2)\bigl)\ =\ 
(1/4)(1 + t y_1 y_2)\ . 
\eqno (3.3) $$
By (3.3) and a simple calculation,
$$ {\rm Corr}\bigl(I(Y_1 = 1), I(Y_2=1)\bigl)\ =\  
\tau(\sigma(Y_1), \sigma(Y_2))\ =\ 
\psi(\sigma(Y_1), \sigma(Y_2))\ =\ t\ . \eqno (3.4) $$
(In working here with the definitions of $\tau(.\, ,.)$ and  
$\psi(.\, ,.)$ in (1.3) and (1.1),
one only needs to check the events of the form
$\{Y_1 = y_1\}$ and $\{Y_2 = y_2\}$ for 
$y_1, y_2, \in \{-1,1\}$, since
$P(A\cap B) - P(A)P(B) = 0$ whenever either $A$ or $B$
is the event $\Omega$ or $\emptyset$.)
\medskip

   By (3.4), and then by 
(3.1), (3.4) (again), and Theorem 4,
$$ t\ =\ \tau\bigl(\sigma(Y_1), \sigma(Y_2)\bigl)\ \leq\  
\tau\Bigl( \sigma(X_1) \vee \sigma(Y_1),\ 
\sigma(X_2) \vee \sigma(Y_2) \Bigl)\ \leq\ t\ , $$
which forces the equality 
$$\tau\Bigl( \sigma(X_1) \vee \sigma(Y_1),\ 
\sigma(X_2) \vee \sigma(Y_2) \Bigl)\ =\ t\ .  
\eqno (3.5)$$
Also, by (3.2),  
$$\rho\Bigl( \sigma(X_1) \vee \sigma(Y_1),\ 
\sigma(X_2) \vee \sigma(Y_2) \Bigl)\ \geq\  
\rho(\sigma(X_1), \sigma(X_2))\ >\  
[t \cdot (1 - \log t)] - \varepsilon\ . \eqno (3.6) $$  
Letting ${\cal A} := \sigma(X_1) \vee \sigma(Y_1)$ and
${\cal B} := \sigma(X_2) \vee \sigma(Y_2)$, one now
obtains Corollary 5 from (3.4), (3.5), and (3.6).
\hfil\break

\noindent {\bf 4. Proof of Theorem 6}
\bigskip

   In the construction for the proof of Theorem 6,
a key role will be played by the example of 
Peyre [20, Theorem 4.1]
described in Theorem 2(II), via the trivially embellished 
form in (the proof of) Corollary 5.
To set that process up, the following
technical lemma (involving just basic calculus)
will be proved first.
\hfil\break

   {\bf Lemma 7.}\ \ {\sl For every $t \in (0,1)$, one 
has that $ t(1 - \log t) > \sin((\pi/2)t)$.}
\hfil\break

   {\bf Proof of Lemma 7.}\ \ With again the usual 
convention $0 \log 0 := 0$, define the function 
$f:[0, \infty) \to {\bf R}$ as follows:  For 
$t \in [0, \infty)$,
$$ f(t)\ :=\ t(1 - \log t)\ -\ \sin((\pi/2)t)\ . 
$$
This function $f$ is continuous on $[0,\infty)$
and has continuous derivatives of all orders on
the open half line $(0,\infty)$.
For $t \in (0,\infty)$, its first three
derivatives are as follows:
$$ \eqalignno{
f'(t)\ &=\ -\log t\ -\ (\pi/2)\cos((\pi/2)t)\ ; \cr
f''(t)\ &=\ -(1/t)\ +\ (\pi/2)^2\sin((\pi/2)t)\ ; \cr
f'''(t)\ &=\ (1/t^2)\ +\ (\pi/2)^3\cos((\pi/2)t)\ . \cr   
}$$

   Now $f'''(t) > 0$ for every $t \in (0,1]$.
Hence $f''(t)$ is strictly increasing for $t \in (0,1]$.
Also $\lim_{t \to 0+} f''(t) = -\infty$ and
$f''(1) = -1 + (\pi/2)^2 > 0$.
Hence there exists a number 
$c \in (0,1)$ such that
$$ f''(t) < 0\ \ {\rm for}\ t \in (0,c); \quad
f''(c) = 0; \quad {\rm and} \quad
f''(t) > 0\ \ {\rm for}\ t \in (c,1]. \eqno (4.1) $$
  
   By (4.1) $f'(t)$ is strictly increasing for
$t \in [c,1]$.   
Also, $f'(1) = 0$.
Hence $f'(t) < 0$ for every $t \in [c,1)$.
Hence $f$ itself is strictly decreasing on $[c,1]$.
Also, $f(1) = 1-1 = 0$.
Hence 
$$ f(t) > 0\ \ {\rm for\ every}\ t \in [c,1). 
\eqno (4.2) $$

   Now $f(0) = 0$ and (by (4.2)) $f(c) > 0$.  
By (4.1), $f$ is ``concave'' ($-f$ is convex) on
the interval $[0,c]$.
It follows that $f(t) > 0$ for every $t \in (0,c)$.
Combining that with (4.2), one has that 
$f(t) > 0$ for all $t \in (0,1)$.  
Consequently, Lemma 7 holds.
\hfil\break

   {\bf Proof of Theorem 6.}\ \  We can (and will) 
let $(\Omega, {\cal F}, P)$ be a probability space 
``rich'' enough to accommodate all 
random variables defined below.
\medskip

   Suppose (1.13) holds.
Applying Lemma 7, let $\varepsilon > 0$ be fixed
sufficiently small that
$$ t(1 - \log t) - \varepsilon > \sin((\pi/2)t)\ . 
\eqno (4.3) $$
(Referring to (1.13), note that both sides of (4.3) are positive.)
\medskip
 
It is well known that (with respect to their
respective Borel $\sigma$-fields), the sets
${\bf R} \times {\bf R}$ and ${\bf R}$ are
bimeasurably isomorphic.
Referring to the final paragraph of Section 3 (i.e.\ 
the final paragraph of the proof of Corollary 5) and
using such an isomorphism, let $(V,W)$ be a random
vector such that
$$ \eqalignno{ 
\tau(\sigma(V), \sigma(W))\ &=\ t 
\quad {\rm and} & (4.4)   \cr
\rho(\sigma(V), \sigma(W))\ &>\ 
\ t(1 - \log t) - \varepsilon\ . & (4.5) \cr  
}$$

Applying (4.5), let $g : {\bf R} \to {\bf R}$ and
$h: {\bf R} \to {\bf R}$ be bounded Borel functions
such that 
$$ r\ :=\ {\rm Corr}(g(V), h(W))\ >\ 
\ t(1 - \log t) - \varepsilon. 
\eqno (4.6) $$
(The existence of such functions $f$ and $g$ with which
one can make ${\rm Corr}(g(V), h(W))$ arbitrarily close
to $\rho(\sigma(V), \sigma(W))$, is a well known
measure-theoretic fact.)
Of course the correlation in (4.6) is positive
(by the sentence after (4.3)). 
We can (and will) normalize those bounded functions 
$g$ and $h$ so that one also has
$$ E[g(V)] = E[h(W)] = 0 \quad {\rm and} \quad
{\rm Var}[g(V)] = {\rm Var}[h(W)] = 1\ . \eqno (4.7) $$  
\medskip

   Let $((V_1, W_1), (V_2, W_2), (V_3, W_3), \dots)$
be a sequence of independent random vectors
such that for each $n \geq 1$, the distribution of
the random vector $(V_n,W_n)$ is the same as that
of $(V,W)$.
\medskip

   For each $n \geq 1$, define the random vector
$(Y_n, Z_n)$ as follows:
$$ Y_n\ :=\ n^{-1/2} \sum_{k=1}^n g(V_k) 
\quad {\rm and} \quad
Z_n\ :=\ n^{-1/2} \sum_{k=1}^n h(W_k)\ .  \eqno (4.8) $$ 
By (4.6), (4.7), (4.8), and the classic central
limit theorem for independent, identically distributed
random vectors whose coordinates have finite second 
moments, one has that
$$ (Y_n, Z_n)\ \Longrightarrow\ (Y,Z) \quad 
{\rm as}\ n \to \infty,     \eqno (4.9) $$
where (i) the symbol $\Longrightarrow$ denotes 
convergence in distribution, and 
(ii) $(Y,Z)$ is a bivariate normal
random vector such that $EY = EZ = 0$,
${\rm Var}(Y) = {\rm Var}(Z) = 1$, and
${\rm Corr}(Y,Z) = r$ (where $r$ is as in (4.6)).
In particular,  
$$ \eqalignno{ 
P(Y_n > 0) &\to 1/2 \quad {\rm as}\ n \to \infty, 
& (4.10) \cr
P(Z_n > 0) &\to 1/2 \quad {\rm as}\ n \to \infty,
\quad {\rm and} 
& (4.11) \cr
P(\{Y_n > 0\} \cap \{Z_n > 0\}) & \to 
P(\{Y > 0\} \cap \{Z > 0\}) \quad {\rm as}\ n \to \infty.
& (4.12) \cr
}$$
By a well known standard calculation involving bivariate
normal distributions, one has that
$$P(\{Y > 0\} \cap \{Z > 0\}) = (1/4) + 
(2\pi)^{-1}\arcsin r\ . $$
(See e.g.\ 
Bradley [2007, Theorem A902 in the Appendix].)\ \
Hence by (4.10), (4.11), and (4.12),  
$$\lim_{n \to \infty}
[P(\{Y_n > 0\} \cap \{Z_n > 0\}) - 
P(Y_n > 0)P(Z_n > 0)]\ = \ (2\pi)^{-1}\arcsin r\ . $$
Hence by (4.10), (4.11), and a simple calculation,
$$ \lim_{n \to \infty} {\rm Corr}(I(Y_n>0), I(Z_n>0))\    
=\ (2/\pi) \arcsin r\ . $$
Hence by (4.8), (4.6), and (4.3) (and (1.13)), one has 
that for all sufficiently large positive integers $n$,
$$ \eqalignno{ 
\tau(\sigma(V_1,V_2, &\dots, V_n),
\sigma(W_1,W_2,\dots,W_n))\ 
\geq\ 
{\rm Corr}(I(Y_n>0), I(Z_n>0)) \cr 
&>\ 
(2/\pi) \arcsin [t(1 - \log t) - \varepsilon] \cr
&>\ (2/\pi) \arcsin [\sin((\pi/2)t)]\    
=\ (2/\pi) \cdot (\pi/2)t\ =\ t\ . & (4.13)  
\cr }$$  

   Now 
$\tau(\sigma(V_1),\sigma(W_1)) 
= \tau(\sigma(V),\sigma(W)) = t$ 
by (4.4).
(Hence for any positive integer $n$, the very
first term in (4.13) is trivially bounded below by
$t$.)\ \ 
Referring to the entire sentence containing (4.13),
let $m$ be the greatest positive integer such that
$\tau(\sigma(V_1,V_2, \dots, V_m),
\sigma(W_1, W_2, \dots, W_m)) = t$.
Define the $\sigma$-fields
${\cal A}_1 := \sigma(V_1,V_2, \dots, V_m)$,
${\cal B}_1 := \sigma(W_1,W_2, \dots, W_m)$,
${\cal A}_2 := \sigma(V_{m+1})$, and
${\cal B}_2 := \sigma(W_{m+1})$.
Then (again see (4.4)), eqs.\ (1.14), (1.15), and (1.16) 
hold.  
That completes the proof of Theorem 6.
\hfil\break

\centerline {REFERENCES}
\bigskip

\def\refs{\medskip\hangindent=25pt\hangafter=1\noindent}

\refs
[1] J.R.\ Blum, D.L.\ Hanson, and L.H.\ Koopmans.  
On the strong law of large numbers for a class of 
stochastic processes. 
{\it Z.\ Wahrsch.\ verw.\ Gebiete\/}\ 2 (1963) 1-11. 

\refs
[2] R.C.\ Bradley.  
Ph.D.\ Dissertation.
University of California at San Diego, 
La Jolla, California, 1978.

\refs
[3] R.C.\ Bradley.  
Equivalent measures of dependence.  
{\it J.\ Multivariate\ Anal.}\ 13 (1983) 167-176.

\refs
[4] R.C.\ Bradley.
On $\rho$-mixing except on small sets.
{\it Pacific J.\ Math.}\ 146 (1990) 217-226.

\refs
[5] R.C.\ Bradley.
Every ``lower psi-mixing'' Markov chain is
``interlaced rho-mixing''.
{\it Stochastic Process.\ Appl.} 72 (1997) 221-239.

\refs
[6] R.C.\ Bradley.
{\it Introduction to Strong Mixing Conditions\/}, 
Vols.\ 1, 2, and 3.
Kendrick Press, Heber City (Utah), 2007.
 
\refs
[7] R.C.\ Bradley.
On mixing properties of some INAR models.
{\it Zapiski Nauchnyh Seminarov POMI\/} 441 (2015) 56-72.
 
\refs
[8] R.C.\ Bradley.
On mixing properties of reversible Markov chains.
{\it New Zealand J.\ Math.\/} 
45 (2015) 71-87.

\refs
[9] R.C.\ Bradley and W.\ Bryc.  
Multilinear forms and measures of dependence between 
random variables.
{\it J.\ Multivariate\ Anal.}\ 16 (1985) 335-367.  
   
\refs
[10] R.C.\ Bradley, W.\ Bryc, and S.\ Janson.
On dominations between measures of dependence. 
{\it J.\ Multivariate\ Anal.}\ 23 (1987) 312-329.  
 
\refs
[11] R.C.\ Bradley and M.\ Peligrad.
Invariance principles under a two-part mixing assumption.
{\it Stochastic Process.\ Appl.} 22 (1986) 271-289.
 
\refs 
[12] A.V.\ Bulinskii.  
On mixing conditions of random fields.  
{\it Theor.\ Probab.\ Appl.}\ 30 (1985) 219-220.
 
\refs
[13] P.\ Cs\'aki and J.\ Fischer.  
On the general notion of maximal correlation.  
{\it Magyar Tud.\ Akad.\ Mat.\ Kutato Int.\ Kozl.}\ 
8 (1963) 27-51.  
 
\refs 
[14] W.\ Doeblin.
Remarques sur la th\'eorie m\'etrique fractions continues.
{\it Compositio Math.}\ 7 (1940) 353-371.
  
\refs
[15] H.O.\ Hirschfeld.  
A connection between correlation and contingency.  
{\it Proc.\ Camb.\ Phil.\ Soc.}\ 31 (1935) 520-524.  

\refs
[16] I.A.\ Ibragimov.
Some limit theorems for stationary processes.
{\it Theor.\ Probab.\ Appl.}\ 7 (1962) 349-382.

\refs
[17] M.\ Iosifescu.
Doeblin and the metric theory of continued fractions: 
A functional theoretic solution to Gauss' 1812 problem.
In: {\it Doeblin and Modern Probability, Proceedings of
the Doeblin Conference held November 2-7, 1991 at the
University of Tubigen's Heinrich Fabri Institut,
Blaubeuren, Germany\/}, (H.\ Cohn, ed.), pp.\ 97-110.
Contemporary Mathematics 149.
American Mathematical Society, Providence (Rhode Island),
1993.

\refs
[18] A.N.\ Kolmogorov and Yu.A.\ Rozanov.  
On strong mixing conditions for stationary Gaussian 
processes.
{\it Theor.\ Probab.\ Appl.\/}\ 5 (1960) 204-208.

\refs
[19] M.\ Peligrad.
On the central limit theorem for weakly dependent sequences
with a decomposed strong mixing coefficient.
{\it Stochastic Process.\ Appl.}\ 42 (1992) 181-193.

\refs 
[20] R.\ Peyre.
Sharp equivalence between $\rho$- and $\tau$-mixing
coefficients.
{\it Studia Math.}\ 216 (2013) 245-270.

\refs
[21] W.\ Philipp.  
The central limit problem for
mixing sequences of random variables.  
{\it Z.\ Wahrsch.\ verw.\ Gebiete\/} 12 (1969) 155-171.
 
\refs
[22] M.\ Rosenblatt.
A central limit theorem and a strong mixing condition.
{\it Proc.\ Natl.\ Acad.\ Sci.\ USA\/} 42 (1956) 43-47.
  
\refs
[23] H.S.\ Witsenhausen.  On sequences of 
pairs of dependent random variables.  
{\it SIAM J.\ Appl.\ Math.}\ 28 (1975) 100-113.

\bye